\documentclass[a4paper,11pt]{amsart}
\usepackage{graphicx}
\usepackage{mathptmx}
\usepackage{amsmath}
\usepackage{amssymb}
\usepackage{enumitem}
\usepackage{xcolor}

\newmuskip\pFqmuskip

\newcommand*\pFq[6][8]{%
  \begingroup % only local assignments
  \pFqmuskip=#1mu\relax
  \mathcode`=\string"8000
  \begingroup\lccode`\~=`\,
  \lowercase{\endgroup\let~}\pFqcomma
  F^{#2}_{#3}{\left(\genfrac..{0pt}{}{#4}{#5}\bigg|#6\right)}%
  \endgroup
}
\newcommand{\pFqcomma}{\mskip\pFqmuskip}

\newtheorem{theorem}{Theorem}
\newtheorem{lemma}[theorem]{Lemma}

\newtheorem{remark}[theorem]{Remark}

\begin{document}

\title[Lah-Bell numbers and polynomials]{Lah-Bell numbers and polynomials}

\author{Dae San Kim}
\address{Department of Mathematics, Sogang University, Seoul 121-742, Republic of Korea}
\email{dskim@sogang.ac.kr}

\author{Taekyun  Kim}
\address{Department of Mathematics, Kwangwoon University, Seoul 139-701, Republic of Korea}
\email{tkkim@kw.ac.kr}

\subjclass[2010]{11B73; 11B83; 05A19}
\keywords{Lah-Bell numbers; Lah-Bell polynomials; bivariate Lah-Bell polynomials; degenerate Lah-Bell numbers; degenerate Lah-Bell polynomials}
 
\maketitle

\begin{abstract}
In this paper, we introduce the Lah-Bell numbers and their natural extensions, namely the Lah-Bell polynomials, and derive some basic properties of such numbers and polynomials by using elementary methods. In addition, we consider the degenerate Lah-Bell numbers and polynomials as degenerate versions of the Lah-Bell numbers and polynomials.
\end{abstract}
 
\section{Introduction}
As is well known, the Stirling number of the second kind $S_{2}(n,k)$, $(n \ge k\ge 0)$, is the number of ways to partition a set with $n$ elements into $k$ non-empty subsets. The $n$-th Bell number $B_{n}$, $(n\ge 0)$ is the number of ways to partition a set with $n$ elements into non-empty subsets. Thus, we have 
\begin{equation}
B_{n}=\sum_{k=0}^{n}S_{2}(n,k),\quad(n\ge 0). \label{1}
\end{equation}  
The Bell polynomials $B_n(x)$ are natural extensions of the Bell numbers (see \eqref{6}, \eqref{10}). \par 
The unsigned Lah number $L(n,k)$ counts the number of ways a set of $n$ elements can be partitioned into $k$ nonempty linearly ordered subsets.
In view of the relationship between the Stirling numbers of the second kind and the Bell numbers,  it is very natural and meaningful to define {\it{the $n$-th Lah-Bell number $B_{n}^{L}$, $(n\ge 0)$, as the number of ways a set of $n$ elements can be partitioned into non-empty linearly ordered subsets}}. Thus, we have 
\begin{equation}
B_{n}^{L}=\sum_{k=0}^{n}L(n,k),\quad(n\ge 0). \label{2}
\end{equation} 
The Lah-Bell polynomials $B_n^{L}(x)$ are also defined as natural extensions of the Lah-Bell numbers (see \eqref{23}, Lemma 4) .\par
The aim of this paper is to study the Lah-Bell numbers and polynomials and to derive some of their basic properties by using elementary methods. We also consider the degenerate Lah-Bell numbers and polynomials (see \eqref{42}, \eqref{44}) as degenerate versions of the Lah-Bell numbers and polynomials. \par
In more detail, for the Lah-Bell numbers and polynomials we derive the generating functions, some relations with Bell numbers and polynomials, and Dobinski-like formulas. We show the connections between Lah numbers and Stirling numbers, recurrence relations and derivatives for the Lah-Bell polynomials. We introduce the bivariate Lah-Bell polynomials, and find the generating function of them and their connection with the bivariate Bell polynomials studied in [15]. In addition, as degenerate versions of the Lah-Bell polynomials and numbers, we introduce the degenerate Lah-Bell polynomials and numbers. Then, for the degenerate Lah-Bell polynomials we deduce the generating function, an explicit expression and their relation with the degenerate Bell polynomials.\par

\vspace{0.1in}

We recall that the falling factorial sequence is given by 
\begin{displaymath}
(x)_{0}=1,\quad (x)_{n}=x(x-1)\cdots(x-n+1), \quad(n \ge 1).
\end{displaymath}
Then, we have 
\begin{equation}
	x^{n}=\sum_{k=0}^{n}S_{2}(n,k)(x)_{k},\quad (n\ge 0),\quad(\mathrm{see}\ [1-15]). \label{3}
\end{equation}
From \eqref{3}, we can derive the following equation. 
\begin{equation}
	\frac{1}{k!}\big(e^{t}-1\big)^{k}=\sum_{n=k}^{\infty}S_{2}(n,k)\frac{t^{n}}{n!},\quad(k\ge 0). \label{4}
\end{equation}
By \eqref{1} and \eqref{4}, we get 
\begin{equation}
e^{(e^{t}-1)}=\sum_{n=0}^{\infty}B_{n}\frac{t^{n}}{n!},\quad(\mathrm{see}\ [4]).\label{5} 
\end{equation}
The Bell polynomials are given by
\begin{equation}
e^{x(e^{t}-1)}=\sum_{n=0}^{\infty}B_{n}(x)\frac{t^{n}}{n!},\quad(\mathrm{see}\ [3,7,8]).\label{6}
\end{equation}
When $x=1$, $B_{n}=B_{n}(1)$ are the $n$-th Bell numbers. 
From \eqref{6}, we have 
\begin{equation}
B_{n}(x)=e^{-x}\sum_{k=0}^{\infty}\frac{k^{n}}{k!}x^{k},\quad(n\ge 0),\quad(\mathrm{see}\ [3,9,12]). \label{7}
\end{equation} \par
It is well known that the Stirling number of the first kind $S_{1}(n,k)$ counts the number of permutations of $n$ elements consisting of $k$ disjoint cycles. 
As an inversion formula of \eqref{3}, we have 
\begin{equation}
(x)_{n}=\sum_{k=0}^{n}S_{1}(n,k)x^{k}. \label{8}
\end{equation}
From \eqref{8}, we note that 
\begin{equation}
\frac{1}{k!}\big(\log(1+t)\big)^k=\sum_{n=k}^{\infty}S_{1}(n,k)\frac{t^{n}}{n!},\quad(k\ge 0),\quad(\mathrm{see}\ [6-8,10-12]). \label{9}
\end{equation}
By \eqref{4} and \eqref{6}, we get 
\begin{equation}
B_{n}(x)=\sum_{K=0}^{n}S_{2}(n,k)x^{k},\quad(n\ge 0),\quad(\mathrm{see}\ [4]). \label{10}
\end{equation} \par
For $n,k\ge 0$, the unsigned Lah numbers are given by 
\begin{equation}
L(n,k)=\binom{n-1}{k-1}\frac{n!}{k!},\quad(\mathrm{see}\ [3,10,12,14]). \label{11}
\end{equation}
The rising factorial sequence is given by 
\begin{equation}
\langle x\rangle_{0}=1,\quad \langle x\rangle_{n}=x(x+1)\cdots(x+n-1),\quad(n\ge 1). \label{12} 
\end{equation}
From \eqref{11}, we note that 
\begin{align}
&\ \langle x\rangle_{n}=\sum_{k=0}^{n}L(n,k)(x)_{k},\label{13} \\
&\ (x)_{n}=\sum_{k=0}^{n}(-1)^{n-k}L(n,k)\langle x\rangle_{k}, \label{14} \\
&\ L(n,k)=\sum_{j=n}^{n}(-1)^{n-j}S_{1}(n,j)S_{2}(j,k), \label{15} \\
&\ L(n,k)=\binom{n}{k}\binom{n-1}{k-1}(n-k)!=\bigg(\frac{n!}{k!}\bigg)^{2}\frac{k}{n(n-k)!},\label{16}\\
&\ L(n,k+1)=\frac{n-k}{k(k+1)}L(n,k),\quad(n \ge 0, k \ge 1). \label{17}
\end{align}

\section{Lah-Bell numbers and polynomials}

In this section, we study the Lah-Bell numbers $B_{n}^{L}$ in \eqref{2} and their extension, namely the Lah-Bell polynomials $B_{n}^{L}(x)$ (see \eqref{23}). In addition, we investigate the degenerate Lah-Bell numbers $B_{n,\lambda}^{L}$

By using \eqref{11}, we easily get 
\begin{equation}
\frac{1}{k!}\bigg(\frac{t}{1-t}\bigg)^{k}=\sum_{n=k}^{\infty}L(n,k)\frac{t^{n}}{n!},\quad(k\ge 0). \label{18}
\end{equation}
From \eqref{2} and \eqref{18}, we note that 
\begin{align}
e^{\frac{t}{1-t}}\ &=\ \sum_{k=0}^{\infty}\frac{1}{k!}\bigg(\frac{t}{1-t}\bigg)^{k}\ =\ \sum_{k=0}^{\infty}\sum_{n=k}^{\infty}L(n,k)\frac{t^{n}}{n!}\label{19} \\
&=\ \sum_{n=0}^{\infty}\bigg(\sum_{k=0}^{n}L(n,k)\bigg)\frac{t^{n}}{n!}\ =\ \sum_{n=0}^{\infty}B_{n}^{L}(n,k)\frac{t^{n}}{n!}. \nonumber	
\end{align}
Therefore, by \eqref{19}, we obtain the following lemma. 
\begin{lemma}
The generating function of Lah-Bell numbers is given by 
\begin{displaymath}
e^{\big(\frac{1}{1-t}-1\big)}=\sum_{n=0}^{\infty}B_{n}^{L}\frac{t^{n}}{n!}. 
\end{displaymath}
\end{lemma}
Replacing $t$ by $1-e^{-t}$ in Lemma 1, we get 
\begin{align}
e^{e^{t}-1}\ &=\ \sum_{k=0}^{\infty}B_{k}^{L}\frac{1}{k!}\big(1-e^{-t}\big)^{k}\nonumber \\
&=\ \sum_{k=0}^{\infty}(-1)^{k}B_{k}^{L}\sum_{n=k}^{\infty}S_{2}(n,k)(-1)^{n}\frac{t^{n}}{n!} \label{20} \\
&=\ \sum_{n=0}^{\infty}\bigg(\sum_{k=0}^{n}(-1)^{n-k}B_{k}^{L}S_{2}(n,k)\bigg)\frac{t^{n}}{n!}. \nonumber
\end{align}
On the other hand, 
\begin{equation}
e^{e^{t}-1}=\sum_{n=0}^{\infty}B_{n}\frac{t^{n}}{n!}. \label{21}
\end{equation}
Therefore, by \eqref{20} and \eqref{21}, we obtain the following theorem.  
\begin{theorem}
For $n\ge 0$, we have 
\begin{displaymath}
B_{n}=\sum_{k=0}^{n}(-1)^{n-k}B_{k}^{L}S_{2}(n,k). 
\end{displaymath}
\end{theorem}
We observe that 
\begin{align}
e^{\frac{1}{1-t}-1}\ &=\ \frac{1}{e}\sum_{k=0}^{\infty}\frac{1}{k!}(1-t)^{-k}\label{22}\\ 
&=\ \frac{1}{e}\sum_{k=0}^{\infty}\frac{1}{k!}\sum_{n=0}^{\infty}\langle k\rangle_{n}\frac{t^{n}}{n!} \nonumber \\
&=\ \sum_{n=0}^{\infty}\bigg\{\frac{1}{e}\sum_{k=0}^{\infty}\frac{\langle k\rangle_{n}}{k!}\bigg\}\frac{t^{n}}{n!}. \nonumber
\end{align}
Therefore, by Lemma 1 and \eqref{22}, we obtain the following Dobinski-like formula for Lah-Bell numbers. 
\begin{theorem}
For $n\ge 0$, we have 
\begin{displaymath}
B_{n}^{L}=\frac{1}{e}\sum_{k=0}^{\infty}\frac{\langle k\rangle_{n}}{k!}. 
\end{displaymath}
\end{theorem}
In view of \eqref{10}, we define Lah-Bell polynomials by 
\begin{equation}
B_{n}^{L}(x)=\sum_{l=0}^{n}L(n,l)x^{l},\quad(n\ge 0). \label{23}
\end{equation}
From \eqref{23}, we note that 
\begin{align}
\sum_{n=0}^{\infty}B_{n}^{L}(x)\frac{t^{n}}{n!}\ &=\ \sum_{n=0}^{\infty}\bigg(\sum_{k=0}^{n}x^{k}L(n,k)\bigg)\frac{t^{n}}{n!} \label{24} \\
&=\ \sum_{k=0}^{\infty}x^{k}\bigg(\sum_{n=k}^{\infty}L(n,k)\frac{t^{n}}{n!}\bigg) \nonumber \\
&=\ \sum_{k=0}^{\infty}x^{k}\frac{1}{k!}\bigg(\frac{1}{1-t}-1\bigg)^{k}\nonumber \\
&=\ e^{x\big(\frac{1}{1-t}-1\big)}.\nonumber
\end{align}
This shows the following result.
\begin{lemma}
The generating function of Lah-Bell polynomials is given by 
\begin{displaymath}
e^{x\big(\frac{1}{1-t}-1\big)}=\sum_{n=0}^{\infty}B_{n}^{L}(x)\frac{t^{n}}{n!}. 
\end{displaymath}
\end{lemma}
Replacing $t$ by $1-e^{-t}$ in Lemma 4, we get 
\begin{align}
e^{x(e^{t}-1)}\ &=\ \sum_{k=0}^{\infty}B_{k}^{L}(x)\frac{1}{k!}\big(1-e^{-t}\big)^{k} \nonumber \\
&=\ \sum_{k=0}^{\infty}(-1)^{k}B_{k}^{L}(x)\sum_{n=k}^{\infty}S_{2}(n,k)(-1)^{n}\frac{t^{n}}{n!}\label{25} \\
&=\ \sum_{n=0}^{\infty}\bigg(\sum_{k=0}^{n}(-1)^{n-k}S_{2}(n,k)B_{k}^{L}(x)\bigg)\frac{t^{n}}{n!}. \nonumber
\end{align}
Therefore, by \eqref{6} and \eqref{25}, we obtain the following theorem. 
\begin{theorem}
For $n\ge 0$, we have 
\begin{displaymath}
B_{n}(x)=\sum_{k=0}^{n}(-1)^{n-k}S_{2}(n,k)B_{k}^{L}(x). 
\end{displaymath}
\end{theorem}
By Theorem 3 and Lemma 4, we get 
\begin{align}
e^{x\big(\frac{1}{1-t}-1\big)}\ &=\ e^{-x}\sum_{k=0}^{\infty}\frac{x^{k}}{k!}(1-t)^{-k}\nonumber \\
&=\ \sum_{k=0}^{\infty}\bigg(e^{-x}\sum_{k=0}^{\infty}\frac{\langle k\rangle_{n}}{k!}x^{k}\bigg)\frac{t^{n}}{n!}.\label{26}	
\end{align}
Therefore, by Lemma 4 and \eqref{26}, we obtain the following Doinski-like formula for Lah-Bell polynomials. 
\begin{theorem}
For $n\ge 0$, the following Dobinski-like formula holds: 
\begin{displaymath}
B_{n}^{L}(x)=e^{-x}\sum_{k=0}^{\infty}\frac{\langle k\rangle_{n}}{k!}x^{k}. 
\end{displaymath}
\end{theorem}
Replacing $t$ by $-\log(1-t)$ in \eqref{6}, we get 
\begin{align}
e^{x\big(\frac{1}{1-t}-1\big)}\ &=\ \sum_{k=0}^{\infty}B_{k}(x)\frac{1}{k!}\big(-\log(1-t)\big)^{k} \nonumber \\
&=\ \sum_{k=0}^{\infty}(-1)^{k}B_{k}(x)\sum_{n=k}^{\infty}(-1)^{n}S_{1}(n,k)\frac{t^{n}}{n!}\label{27} \\
&=\ \sum_{n=0}^{\infty}\bigg(\sum_{k=0}^{n}(-1)^{n-k}S_{1}(n,k)B_{k}(x)\bigg)\frac{t^{n}}{n!}.\nonumber
\end{align}
Therefore, by Lemma 4 and \eqref{27}, we obtain the following theorem. 
\begin{theorem}
For $n\ge 0$, we have 
\begin{displaymath}
B_{n}^{L}(x)=\sum_{k=0}^{n}(-1)^{n-k}S_{1}(n,k)B_{k}(x). 
\end{displaymath}
\end{theorem}
Replacing $t$ by $-\log(1-t)$ in \eqref{4}, we have 
\begin{align}
\frac{1}{k!}\bigg(\frac{1}{1-t}-1\bigg)^{k}\ &=\ \sum_{l=k}^{\infty}S_{2}(l,k)\frac{1}{l!}\big(-\log(1-t)\big)^{l}\label{28} \\
&=\ \sum_{l=k}^{\infty}(-1)^{l}S_{2}(l,k)\sum_{n=l}^{\infty}S_{1}(n,l)(-1)^{n}\frac{t^{n}}{n!} \nonumber \\
&=\ \sum_{n=k}^{\infty}\bigg(\sum_{l=k}^{n}(-1)^{n-l}S_{1}(n,l)S_{2}(l,k)\bigg)\frac{t^{n}}{n!}\nonumber.
\end{align}
On the other hand, 
\begin{equation}
\frac{1}{k!}\bigg(\frac{1}{1-t}-1\bigg)^{k}=\sum_{n=k}^{\infty}L(n,k)\frac{t^{n}}{n!}. \label{29}
\end{equation}
By \eqref{28} and \eqref{29}, we get 
\begin{equation}
L(n,k)=\sum_{l=k}^{n}(-1)^{n-l}S_{1}(n,l)S_{2}(l,k),\quad(n\ge 0). \label{30}
\end{equation}
By replacing $t$ by $1-e^{-t}$ in \eqref{29}, we get 
\begin{align}
\frac{1}{k!}\big(e^{t}-1\big)^{k}\ &=\ \sum_{l=k}^{\infty}L(l,k)\frac{1}{l!}\big(1-e^{-t}\big)^{l}\label{31}\\
&=\ \sum_{l=k}^{\infty}(-1)^{l}L(l,k)\sum_{n=l}^{\infty}S_{2}(n,l)(-1)^{n}\frac{t^{n}}{n!} \nonumber \\
&=\ \sum_{n=k}^{\infty}\bigg(\sum_{l=k}^{n}(-1)^{n-l}S_{2}(n,l)L(l,k)\bigg)\frac{t^{n}}{n!}.\nonumber
\end{align}
Therefore, by \eqref{4} and \eqref{31}, we get 
\begin{equation}
S_{2}(n,k)=\sum_{l=k}^{n}(-1)^{n-l}S_{2}(n,l)L(l,k), 
\end{equation}
where $n,k\ge 0$, with $n\ge k$. 
\begin{theorem}
For $n,k\ge 0$, with $n\ge k$, we have 
\begin{displaymath}
S_{2}(n,k)=\sum_{l=k}^{n}(-1)^{n-l}S_{2}(n,l)L(l,k),
\end{displaymath} 
and 
\begin{displaymath}
L(n,k)=\sum_{l=k}^{n}(-1)^{n-l}S_{1}(n,l)S_{2}(l,k).
\end{displaymath}
\end{theorem}
From Lemma 4, we note that 
\begin{equation}
x\frac{d}{dt}\bigg(\frac{1}{1-t}\bigg)e^{x\big(\frac{1}{1-t}-1\big)}
=\ \sum_{n=0}^{\infty}B_{n+1}^{L}(x)\frac{t^{n}}{n!}. \label{33}
\end{equation}
On the other hand,
\begin{align}
x\frac{d}{dt}\bigg(\frac{1}{1-t}\bigg)e^{x\big(\frac{1}{1-t}-1\big)}\ &=\ x\sum_{l=0}^{\infty}(l+1)!\frac{t^{l}}{l!}\sum_{m=0}^{\infty}B_{m}^{L}(x)\frac{t^{m}}{m!} \label{34} \\
&=\ \sum_{n=0}^{\infty}\bigg(x\sum_{m=0}^{n}\binom{n}{m}(n-m+1)!B_{m}^{L}(x)\bigg)\frac{t^{n}}{n!}. \nonumber
\end{align}
Therefore, by \eqref{33} and \eqref{34}, we obtain the following theorem. 
\begin{theorem}
For $n\ge 0$, we have 
\begin{displaymath}
B_{n+1}^{L}(x)=x\sum_{m=0}^{n}\binom{n}{m}(n-m+1)!B_{m}^{L}(x).	
\end{displaymath}
\end{theorem}
Again, from Lemma 4, we observe that 
\begin{align*}
\sum_{n=1}^{\infty}\frac{d}{dx}B_{n}^{L}(x)\frac{t^{n}}{n!}\ &=\ \frac{d}{dx}e^{x\big(\frac{1}{1-t}-1\big)}\ =\ \bigg(\frac{1}{1-t}-1\bigg)e^{x\big(\frac{1}{1-t}-1\big)}\\
&=\ \sum_{l=1}^{\infty} l! \frac{t^{l}}{l!}\sum_{m=0}^{\infty}B_{m}^{L}(x)\frac{t^{m}}{m!}\\
&=\ \sum_{n=1}^{\infty}\bigg(\sum_{m=0}^{n-1}\binom{n}{m}(n-m)!B_{m}^{L}(x)\bigg)\frac{t^{n}}{n!}.  
\end{align*}
Thus, we have 
\begin{equation}
\frac{d}{dx}B_{n}^{L}(x)=\sum_{m=0}^{n-1}\binom{n}{m}(n-m)!B_{m}^{L}(x),\quad(n\ge 1). \label{35}
\end{equation}
Therefore, by \eqref{35}, we obtain the following theorem.
\begin{theorem}
For $n\ge 1$, we have 
\begin{displaymath}
\frac{d}{dx}B_{n}^{L}(x)=\sum_{m=0}^{n-1}\binom{n}{m}(n-m)!B_{m}^{L}(x).
\end{displaymath}
\end{theorem}
For $n\ge 0$, the bivariate Bell polynomials are defined by 
\begin{equation}
B_{n}(x,y)=\sum_{k=0}^{n}S_{2}(n,k)(x)_{k}y^{k},\quad(\mathrm{see}\ [15]). \label{36}
\end{equation}
Letting $y\rightarrow y/x$ and then $x\rightarrow\infty$, we see that the bivariate Bell polynomial $B_{n}(x,y)$ reduces to the univariate Bell polynomial $B_{n}(y)$.\par 
From \eqref{36}, we note that 
\begin{equation}
\sum_{n=0}^{\infty}B_{n}(x,y)\frac{t^{n}}{n!}\ =\ \big(1+y(e^{t}-1)\big)^{x},\quad(\mathrm{see}\ [15]). \label{37}
\end{equation}
In view of \eqref{36}, we define bivariate Lah-Bell polynomials by
\begin{equation}
B_{n}^{L}(x,y)=\sum_{k=0}^{n}L(n,k)(x)_{k}y^{k},\quad(n\ge 0). \label{38}
\end{equation} 
From, \eqref{38}, we note that 
\begin{align}
\sum_{n=0}^{\infty}B_{n}^{L}(x,y)\frac{t^{n}}{n!}\ &=\ \sum_{n=0}^{\infty}\bigg(\sum_{k=0}^{n}L(n,k)(x)_{k}y^{k}\bigg)\frac{t^{n}}{n!} \label{39}\\
&=\ \sum_{k=0}^{\infty}(x)_{k}y^{k}\sum_{n=k}^{\infty}	L(n,k)\frac{t^{n}}{n!} \nonumber \\
&=\ \sum_{k=0}^{\infty}\binom{x}{k}y^{k}\bigg(\frac{1}{1-t}-1\bigg)^{k}\nonumber \\
&=\ \bigg(1+y\bigg(\frac{1}{1-t}-1\bigg)\bigg)^{x}. \nonumber
\end{align}
Therefore, by \eqref{39}, we obtain the following lemma. 
\begin{lemma}
The generating function of bivariate Lah-Bell polynomais is given by 
\begin{displaymath}
\bigg(1+y\bigg(\frac{1}{1-t}-1\bigg)\bigg)^{x}=\sum_{n=0}^{\infty}B_{n}^{L}(x,y)\frac{t^{n}}{n!}. 
\end{displaymath}
\end{lemma}

\begin{remark} Letting $y\rightarrow y/x$ and $x\rightarrow\infty$, we note the bivariate Lah-Bell polynomial $B_{n}^{L}(x,y)$ reduces to the univariate Lah-Bell polynomial $B_{n}^{L}(y)$, $(n\ge 0)$.  
\end{remark}
Replacing $t$ by $1-e^{-t}$ in Lemma 11, we get
\begin{align}
\big(1+y(e^{t}-1)\big)^{x}\ &=\ \sum_{k=0}^{\infty}B_{k}^{L}(x,y)\frac{1}{k!}\big(1-e^{-t}\big)^{k} \nonumber \\
&=\ \sum_{k=0}^{\infty}(-1)^{k}B_{k}^{L}(x,y)\sum_{n=k}^{\infty}S_{2}(n,k)\frac{(-t)^{n}}{n!}\label{40}\\
&=\ \sum_{n=0}^{\infty}\bigg(\sum_{k=0}^{n}(-1)^{n-k}S_{2}(n,k)B_{k}^{L}(x,y)\bigg)\frac{t^{n}}{n!}. \nonumber
\end{align}
Replacing $t$ by $-\log(1-t)$ in \eqref{37}, we get 
\begin{align}
\bigg(1+y\bigg(\frac{1}{1-t}-1\bigg)\bigg)^{x}\ &=\ \sum_{k=0}^{\infty}B_{k}(x,y)(-1)^{k}\frac{1}{k!}\big(\log(1-t)\big)^{k} \label{41} \\
&=\ \sum_{k=0}^{\infty}(-1)^{k}B_{k}(x,y)\sum_{n=k}^{\infty}(-1)^{n}S_{1}(n,k)\frac{t^{n}}{n!}\nonumber \\
&=\ \sum_{n=0}^{\infty}\bigg(\sum_{k=0}^{n}(-1)^{n-k}S_{1}(n,k)B_{k}(x,y)\bigg)\frac{t^{n}}{n!}\nonumber.\end{align}
Therefore, by \eqref{37}, Lemma 11, \eqref{40} and \eqref{41}, we obtain the following theorem. 
\begin{theorem}
For $n\ge 0$, we have 
\begin{displaymath}
B_{n}^{L}(x,y)=\sum_{k=0}^{n}(-1)^{n-k}S_{1}(n,k)B_{k}(x,y), 
\end{displaymath}
and 
\begin{displaymath}
B_{n}(x,y)=\sum_{k=0}^{n}(-1)^{n-k}S_{2}(n,k)B_{k}^{L}(x,y). 
\end{displaymath}
\end{theorem}
For any $0 \ne \lambda\in\mathbb{R}$, the degenerate exponential functions are defined by 
\begin{displaymath}
e_{\lambda}^{x}(t)=(1+\lambda t)^{\frac{x}{\lambda}}=\sum_{n=0}^{\infty}(x)_{n,\lambda}\frac{t^{n}}{n!},\quad(\mathrm{see}\ [\ 7,9]),
\end{displaymath}
where $(x)_{0,\lambda}=1,\ (x)_{n,\lambda}=x(x-\lambda)\cdots(x-(n-1)\lambda),\ (n\ge 1)$. \\
When $x=1$, $e_{\lambda}(t)=e_{\lambda}^{1}(t)=\sum_{n=0}^{\infty}\frac{(1)_{n,\lambda}}{n!}t^{n}$.
Note that $\lim_{\lambda\rightarrow 0}e_{\lambda}(t)=e^{t}$. \par
Now, we define the degenerate Lah-Bell polynomials by
\begin{equation}
e_{\lambda}^{x}\bigg(\frac{1}{1-t}-1\bigg)=\sum_{n=0}^{\infty}B_{n,\lambda}^{L}(x)\frac{t^{n}}{n!}. \label{42}
\end{equation}
When $x=1$, $B_{n,\lambda}^{L}=B_{n,\lambda}^{L}(1)$ are called the degenerate Lah-Bell numbers. \par 
From \eqref{42}, we note that 
\begin{align}
\sum_{n=0}^{\infty}B_{n,\lambda}^{L}(x)\frac{t^{n}}{n!}\ &=\ \sum_{k=0}^{\infty}(x)_{k,\lambda}\frac{1}{k!}\bigg(\frac{1}{1-t}-1\bigg)^{k}\nonumber  \\
&=\ \sum_{k=0}^{\infty}(x)_{k,\lambda}\sum_{n=k}^{\infty}L(n,k)\frac{t^{n}}{n!} \label{43} \\
&=\ \sum_{n=0}^{\infty}\bigg(\sum_{k=0}^{n}L(n,k)(x)_{k,\lambda}\bigg)\frac{t^{n}}{n!}. \nonumber
\end{align}
By comparing the coefficients on both sides of \eqref{43}, we get 
\begin{equation}
B_{n,\lambda}^{L}(x)=\sum_{k=0}^{n}L(n,k)(x)_{k,\lambda},\quad(n\ge 0). \label{44}
\end{equation}
It is known that the degenerate Bell polynomials are given by
\begin{equation}
e_{\lambda}^{x}\big(e^{t}-1\big)=\sum_{n=0}^{\infty}B_{n,\lambda}(x)\frac{t^{n}}{n!},\quad(\mathrm{see}\ [\ 9]).\label{45}
\end{equation}
Replacing $t$ by $1-e^{-t}$ in \eqref{42}, we get 
\begin{align}
e_{\lambda}^{x}\big(e^{t}-1\big)\ &=\ \sum_{k=0}^{\infty}B_{k,\lambda}^{L}(x)\frac{1}{k!}\big(1-e^{-t}\big)^k \label{46} \\
&=\ \sum_{n=0}^{\infty}\bigg(\sum_{k=0}^{n}(-1)^{n-k}S_{2}(n,k)B_{k,\lambda}^{L}(x)\bigg)\frac{t^{n}}{n!}.\nonumber
\end{align}
From \eqref{45} and \eqref{46}, we note that 
\begin{equation}
B_{n,\lambda}(x)=\sum_{k=0}^{n}(-1)^{n-k}S_{2}(n,k)B_{k,\lambda}^{L}(x),\quad(n\ge 0). \label{47}
\end{equation}
Replacing $t$ by $-\log(1-t)$ in \eqref{45}, we have 
\begin{align}
e_{\lambda}^{x}\bigg(\frac{t}{1-t}-1\bigg)\ &=\ \sum_{k=0}^{\infty}(-1)^{k}B_{k,\lambda}(x)\frac{1}{k!}\big(\log(1-t)\big)^{k} \label{48}  \\
&=\ \sum_{n=0}^{\infty}\bigg(\sum_{k=0}^{n}(-1)^{n-k}S_{1}(n,k)B_{k,\lambda}(x)\bigg)\frac{t^{n}}{n!}.\nonumber
\end{align}
Thus, we have 
\begin{displaymath}
B_{n,\lambda}^{L}(x)=\sum_{k=0}^{n}(-1)^{n-k}S_{1}(n,k)B_{k,\lambda}(x),\quad (n\ge 0). 
\end{displaymath}
\indent Recall that the Laguerre polynomials $L_{n}^{(\alpha)}(x)$ of order $\alpha,\,\,(\alpha >-1)$, are given by (see [12])
\begin{equation}
(1-t)^{-\alpha-1}e^{x\frac{t}{t-1}}=\sum_{n=0}^{\infty}L_{n}^{(\alpha)}(x)\frac{t^n}{n!}.\label{49}
\end{equation}
By Lemma 4 and \eqref{49},
\begin{align*}
(1-t)^{-\alpha-1}&=\sum_{m=0}^{\infty}B_m^{L}(x)\frac{t^m}{m!}\sum_{l=0}^{\infty}L_{l}^{(\alpha)}(x)\frac{t^l}{l!}  \\
&=\sum_{n=0}^{\infty}\bigg(\sum_{m=0}^n\binom{n}{m}B_m^{L}(x)L_{n-m}^{(\alpha)}(x)\bigg)\frac{t^n}{n!},
\end{align*}
which shows that
\begin{equation*}
\langle \alpha+1 \rangle_n=\sum_{m=0}^n\binom{n}{m}B_m^{L}(x)L_{n-m}^{(\alpha)}(x).
\end{equation*}

\section{conclusion}
Taking into account the relationship between the Stirling numbers of the second kind and the Bell numbers and in light of the combinatorial meaning of the unsigned Lah numbers, we introduced the Lah-Bell numbers and their natural extensions, namely the Lah-Bell polynomials. We derived some basic properties of such numbers and polynomials by using elementary methods. We also considered the degenerate Lah-Bell numbers and polynomials as degenerate versions of the Lah-Bell numbers and polynomials. \par
In more detail, for the Lah-Bell numbers and polynomials we derived the generating functions, some relations with Bell numbers and polynomials, and Dobinski-like formulas. We showed the connections between Lah numbers and Stirling numbers, recurrence relations and derivatives for the Lah-Bell polynomials. We introduced the bivariate Lah-Bell polynomials, and found the generating function of them and their connection with the bivariate Bell polynomials. In addition, as degenerate versions of the Lah-Bell polynomials and numbers, we introduced the degenerate Lah-Bell polynomials and numbers. Then, for the degenerate Lah-Bell polynomials we deduced the generating function, an explicit expression and their relation with the degenerate Bell polynomials. \par
In recent years, one of our research areas of study has been to explore some special numbers and polynomials and their degenerate versions, and to discover their arithmetical and combinatorial properties and some of their applications. We would like to continue to work on these by exploiting various means like generating functions, combinatorial methods, $p$-adic analysis, umbral calculus, differential equations and probability theory.

\end{document}